\documentclass{icm2010}

\usepackage{amsfonts}

\ifx\pdfoutput\undefined
  \usepackage[dvips]{graphicx}
   \else
   \usepackage[pdftex]{graphicx}
   \fi
   \usepackage{epstopdf}

\newcommand{\F}{\mathsf F}
\newcommand{\A}{{\mathsf A}}
\newcommand{\B}{{\mathsf B}}
\newcommand{\M}{{\mathsf M}}
\newcommand{\CC}{{\mathbb C}}
\newcommand{\RR}{{\mathbb R}}
\newcommand{\ZZ}{{\mathbb Z}}
\newcommand{\PP}{{\mathbb P}}
\newcommand{\dV}{{\mathsf V}}
\newcommand{\C}{{\mathsf C}}
\newcommand{\eps}{\epsilon}
\newcommand{\dW}{{\mathsf W}}
\newcommand{\X}{{\mathbb X}}
\newcommand{\Y}{{\mathbb Y}}
\newcommand{\Z}{{\mathbb Z}}
\newcommand{\cC}{{\mathcal C}}
\newcommand{\T}{{\mathbb T}}
\newcommand{\cD}{{\mathcal D}}
\newcommand{\Dtwo}{{D_{\ZZ_2}}}
\newcommand{\cM}{{\mathcal M}}
\newcommand{\LG}{{{}^LG}}
\newcommand{\frg}{{\mathfrak g}}
\newcommand{\cO}{{\mathcal O}}

\newcommand{\ra}{\rightarrow}

\newcommand{\Cob}{{\rm Bord}}
\newcommand{\Vect}{{\rm Vect}}
\newcommand{\Hom}{{\rm Hom}}
\newcommand{\Sym}{{\rm Sym}}
\newcommand{\bpartial}{{\bar\partial}}
\newcommand{\bul}{\bullet}

\newtheorem*{conjecture}{Conjecture}

\title[TFT, Higher Categories and Applications]{Topological Field Theory, Higher Categories, and Their Applications}

\author[A. Kapustin]{Anton Kapustin}

\contact[kapustin@theory.caltech.edu]{California Institute of Technology, Pasadena, CA 91125, United States}

\begin{document}

\begin{abstract}

It has been common wisdom among mathematicians that Extended Topological Field Theory in dimensions higher than two is naturally formulated in terms of $n$-categories with $n> 1$. Recently the physical meaning of these higher categorical structures has been recognized and concrete examples of Extended TFTs  have been constructed. Some of these examples, like the Rozansky-Witten model, are of geometric nature, while others are related to representation theory. I outline two application of higher-dimensional TFTs. One is related to the problem of classifying monoidal deformations of the derived category of coherent sheaves, and the other one is geometric Langlands duality.
\end{abstract}



\section{Introduction}

The notion of functional integral\footnote{The term "functional integral" is synonymous with "path-integral", but is more descriptive, since in Quantum Field Theory one integrates over a space of functions of several variables rather than over a space of paths.}  plays a central role in quantum field theory, but it has defied attempts at a rigorous mathematical formulation, except in some special cases (typically in space-time dimension 2). Topological Field Theory (TFT) provides a useful playground for studying properties of the functional integral in a simplified setting and has been the subject of many works since the pioneering papers by E.~Witten \cite{Witten:DonaldsonTFT, Witten:Topsigma, Witten:Jones}. The physical definition of the functional integral uses an ill-defined measure on the space of field configurations, but one can use the usual mathematical ploy and try to axiomatize properties of the functional integral without making a direct reference to this measure. The first attempt at such an axiomatization was made by M.~Atiyah \cite{Atiyah:TFT}. Atiyah defines a TFT in $n$-dimensions as a functor $\F$ from a certain geometrically defined category $\Cob_n$ to the category of complex vector spaces $\Vect$ (or to the category of $\ZZ_2$-graded complex vector spaces $\Vect_{\ZZ_2}$). The category $\Cob_n$ has as its objects compact oriented $n-1$-manifolds without boundary and has as its morphisms oriented bordisms between such manifolds. $\F$ is supposed to be invariant with respect to diffeomorphisms. The disjoint union gives the category $\Cob_n$ a symmetric monoidal structure whose identity object is the empty $n-1$-manifold. The category $\Vect$ also has a natural symmetric monoidal structure given by the tensor product; the identity object is the field $\CC$. The functor $\F$ is required to be monoidal; in particular, it sends the disjoint union of two $n-1$-manifolds $M_1$ and $M_2$ to the tensor product $V(M_1)\otimes V(M_2)$, and it sends the empty manifold to $\CC$. 

Since we can regard a closed oriented $n$-manifold $N$ as a bordism between $\emptyset$ and $\emptyset$, the functor $\F$ sends any such $N$ to a linear map from $\CC$ to $\CC$, i.e. a complex number $\F(N)$. This number is called the partition function of the TFT on the manifolds $N$. So we see that an $n$-dimensional TFT assigns a number to a closed oriented $n$-manifold and a vector space to a closed oriented $n-1$-manifold. 

It is natural to ask if an $n$-dimensional TFT assigns anything to closed oriented manifolds of lower dimensions. An obvious guess is that it assigns a $\CC$-linear category to an $n-2$-manifold, a $\CC$-linear 2-category to an $n-3$-manifold, etc. The resulting gadget is usually called an {\it Extended Topological Field Theory}.\footnote{It is likely that the language of $(\infty,n)$-categories whose theory is being developed by J. Lurie \cite{Lurie:Goodwillie} is even better suited for TFT applications \cite{Lurie:TFT,FHLT}. Its physical significance remains unclear at the time of writiing.}  Extending the TFT functor to lower-dimensional manifolds is natural if we consider gluing closed manifolds out of manifolds with boundaries. For example, given two oriented $n-1$-dimensional manifolds $N_1$ and $N_2$ and an orientation-reversing diffeomorphism $\partial N_1\ra\partial N_2$, we can glue $N_1$ and $N_2$ along their common boundary and get a closed oriented $n-1$-manifold $N_{12}$. An Extended TFT in $n$-dimensions assigns to $\partial N_1$ a $\CC$-linear category $\F(\partial N_1)$ and assigns to $N_1$ and ${\bar N}_2$ (the orientation-reversal of $N_2$) objects $\F(N_1)$ and $\F({\bar N}_2)$ of this category. The fact that $N_{12}$ can be glued from $N_1$ and $N_2$ means that the vector space $\F(N_{12})$ is the space of morphisms from the object $\F(N_1)$ to the object $\F({\bar N}_2)$. 

While the relevance of higher categories for TFT in higher dimensions has been recognized by experts for some time \cite{Freed1,Freed2, Lawrence,BaezDolan}, an axiomatic definition of an Extended TFT has not been formulated for an obvious reason: the lack of a universally accepted definition of an $n$-category for $n>2$.\footnote{We remind that one distinguishes strict and weak $n$-categories. While the former are easily defined, they almost never occur in practice; to define Extended TFTs one needs weak $n$-categories.} This technical obstacle was compounded by a lack of understanding of the physical meaning of higher categories.

The correct definition of a weak $n$-category being non-obvious, one is forced to go back to the physical roots of the subject. We will first discuss two-dimensional TFTs which have been studied extensively because of their connection with Mirror Symmetry and explain why boundary conditions in a 2d TFT form a category. This  observation is due to M.~Douglas \cite{Douglas}. Then we will move on to three dimensions and explain why boundary conditions in a 3d TFT form a 2-category. Applying these observations to $n$-dimensional TFT we will be able to see from a more physical viewpoint why $n$-dimensional TFT assigns a $\CC$-linear $k-1$-category to a compact oriented $n-k$-manifold. Then we will describe two examples of TFTs in three and four dimensions and their applications to two different mathematical problems: the classification of monoidal deformations of the derived category of coherent sheaves and the Geometric Langlands Program.

\section{Extended Topological Field Theory from a physical viewpoint}

\subsection{Extended TFT in two dimensions}

 Consider a 2d TFT on a compact oriented 2-manifold $\Sigma$ with a nonempty boundary. It turns out necessary to impose some conditions on the values of the fields on $\partial\Sigma$ for the functional integral to be well-defined on the physical level of rigor. Roughly speaking, these conditions must define a Lagrangian submanifold in the space of boundary values of the fields, where the symplectic form arises from the boundary terms in the variation of the action. On the classical level, boundary conditions are needed to make the initial-value problem for the classical equations of motion to be well-posed and to ensure the existence of a symplectic form on the space of solutions. 
 
 Boundary conditions in a 2d TFT (also known as branes) form a $\CC$-linear category. This is the category which the 2d TFT assigns to a point. Morphisms in this category are {\it boundary-changing local operators}. To explain informally what this means, suppose $\Sigma$ is a half-plane $\{(x,y)\in\RR^2\vert x \geq 0\}$, and one imposed a boundary condition $\A$ on the half-line $\{(0,y)\vert y<0\}$ and a boundary condition $\B$ on the half-line $\{(0,y)\vert y>0\}$ (see Fig. 1). At the special point $(0,0)$ additional data are needed to specify the functional integral uniquely. These data define a boundary-changing point operator $\cO_{AB}$ between $\A$ and $\B$.

\begin{figure}[htbp] \label{fig1}
\centering
\includegraphics[bb=0 10 350 180,height=2in, scale=2.5, clip=true]{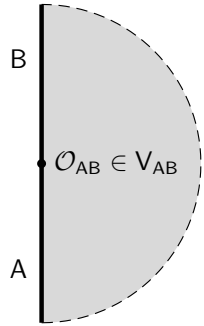}
\caption{Morphisms in the category of boundary conditions are boundary-changing point operators.}
\end{figure}

It is a basic physical principle that the set of boundary-changing point operators has the structure of a (graded) vector space. To see why, let us introduce polar coordinates $(r,\phi)$, $r\in\RR_+$, $\phi\in [-\pi/2,\pi/2]$,  so that the origin $(0,0)$ is given by $r=0$. To avoid dealing with divergences ubiquitous in quantum field theory one may cut out a small half-disc  $r < \eps$ for some $\eps>0$ and replace the boundary-changing point operator by a suitable boundary condition on a semi-circle $r=\eps$ (Fig. 2). (Unlike the boundary conditions $\A$ and $\B$ the boundary condition corresponding to a boundary-changing point operator is nonlocal, in general, in the sense that it does not merely constrain the values of the fields and a finite number of their derivatives along the boundary but may involve constraints on the Fourier components of the restrictions of the fields to the boundary.) The key remark is that the half-plane with a half-disc removed is diffeomorphic to a product of a half-line parameterized by $r$, $r\geq\eps,$ and the interval $[-\pi/2,\pi/2]$ parameterized by $\phi$. We now reinterpret $r$ as the time coordinate and $\phi$ as the spatial coordinate. On the spatial boundaries $\phi=\mp \pi/2$ we have boundary conditions $\A$ and $\B$, while the boundary condition at $r=\eps$ is now regarded as an initial condition. Initial states in any quantum theory form a vector space (in fact, a Hilbert space). We conclude that boundary-changing point operators between boundary conditions $\A$ and $\B$ form a vector space $\dV_{\A\B}$.\footnote{Note that orientation is important here, so $\dV_{\A\B}$ is not the same as $\dV_{\B\A}$.}

\begin{figure}[htbp]  \label{fig2}
\centering
\includegraphics[bb=0 10 350 180,height=2in, scale=2.5, clip=true]{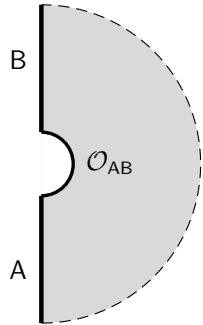}
\caption{A boundary-changing point operator is equivalent to a (possibly nonlocal) boundary condition.}
\end{figure}

Next we consider a situation where $\Sigma$ is a half-plane, but its boundary is divided into three pieces: a half-line $y<0$, an interval $0< y< a$, and a half-line $y>a$. We impose boundary conditions $\A,\B,\C$ on the three pieces respectively, so we need two boundary-changing point operators which are elements of vector spaces $\dV_{\A\B}$ and $\dV_{\B\C}$ (Fig. 3). In a TFT the limit $a\ra 0$ always exists, and one should be able to interpret these two boundary-changing point operators as a single boundary-changing point operator between $\A$ and $\C$, i.e. an element of the vector space $\dV_{\A\C}$. This gives a `fusion product''
$$
\dV_{\A\B}\times \dV_{\B\C}\ra {\dV}_{\A\C}
$$
One may further argue that this product is bilinear and associative in an obvious sense. Note also that for any $\A$ the vector space $\dV_{\A\A}$ has a special element: the boundary-changing point operators which is trivial. This element serves as a unit in the  algebra $\dV_{\A\A}$. Altogether we obtain a category whose objects are branes, whose morphisms are elements of vector spaces $\dV_{\A\B}$, and with the composition of morphisms defined by means of the fusion product.

\begin{figure}[htbp] \label{fig3}
\centering
\includegraphics[bb=100 10 450 200,height=2in, scale=2.5, clip=true]{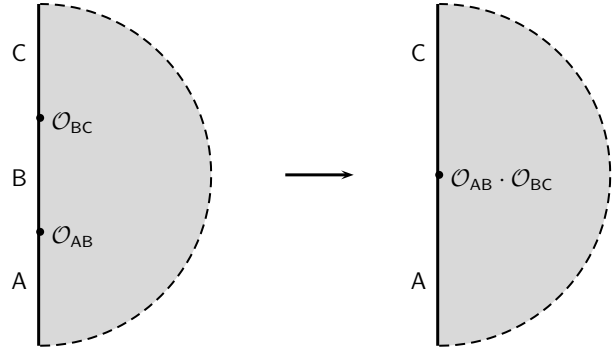}
\caption{Composition of morphisms corresponds to fusing boundary-changing point operators.  Fusion product is denoted by  a dot.}
\end{figure}

Axioms of 2d TFT as usually formulated further imply that the resulting category is self-dual, in the sense that the spaces $\dV_{\A\B}$ and $\dV_{\B\A}$ are naturally dual, but we will not emphasize this aspect, since in some cases with infinite-dimensional spaces $\dV_{\A\B}$ this requirement is not satisfied. 

The last remark we want to make about Extended TFTs in two dimensions is that  instead of closed oriented manifolds one may consider compact oriented manifolds with a nonempty boundary. The connected components of the boundary should be labeled by objects of the category of branes. We will call this labeling a decoration of a manifold. Extended TFT assigns a complex number to a decorated 2-manifold (the value of the functional integral with the corresponding boundary conditions). It assigns a vector space to 
a decorated 1-manifold. The only connected decorated 1-manifold is an interval, whose decoration consists of an ordered pair $(\A,\B)$ of branes. The vector space assigned to such a pair is the space of boundary-changing local operators $\dV_{\A\B}$ introduced above. A good way to think about this rule is the following: if we consider our 2d TFT on a 2-manifold of the form $\RR\times [0,1]$, where the two connected components of the boundary are labeled by $\A$ and $\B$, then we may regard it as a 1d TFT on $\RR$ (this is called Kaluza-Klein reduction).  A 1d TFT is simply a quantum mechanical system, and its Hilbert space of states is what we assign to the interval $[0,1]$. 

\subsection{Extended TFT in three dimensions}

Three-dimensional TFT is supposed to assigns a $\CC$-linear category to a closed oriented 1-manifold and a $\CC$-linear 2-category to a point. Let us first explain the physical meaning of the former. Consider a 3d TFT on a manifold of the form $S^1\times \Sigma$, where $\Sigma$ is an oriented 2-manifold which may be noncompact or have a nonempty boundary. Another basic physical principle (Kaluza-Klein reduction) is that in such a case one can describe the physics of the compactified theory by an {\it effective 2d TFT} on $\Sigma$. By definition, the category assigned to a circle is the category of branes in this effective 2d TFT.

The 2-category assigned to a point is the 2-category of boundary conditions in the 3d TFT. To explain where the 2-category structure comes from, consider a 3d TFT on an oriented 3-manifold $W$ with a nonempty boundary, imagine that a connected component of $\partial W$ is subdivided by closed curves into domains, and that one imposed unrelated boundary conditions on different domains. Each domain is thus labeled by an element of the set of boundary conditions. A closed curve separating the domains labeled by boundary conditions $\X$ and $\Y$ is itself labeled by an element of a set $\dW_{\X\Y}$ which determines how fields behave in the neighborhood of the closed curve. Elements of the set $\dW_{\X\Y}$ will be called boundary-changing line operators from $\X$ to $\Y$. Boundary-changing line operators may be fused together (Fig. 4) which gives rise to a fusion product 
\begin{equation}\label{fusion}
\dW_{\X\Y}\times \dW_{\Y\Z}\ra \dW_{\X\Z},\quad (\A,\B)\mapsto \A\otimes\B,\ \forall \A\in\dW_{\X\Y},\forall \B\in\dW_{\Y\Z}.
\end{equation}

\begin{figure}[htbp]  \label{fig4}
\centering
\includegraphics[bb=100 10 550 200,height=2in, scale=2.5, clip=true]{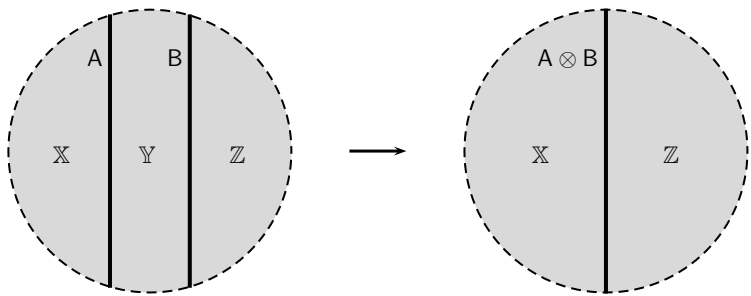}
\caption{Boundary-changing line operator $\A$ from a boundary condition $\X$ to a boundary condition $\Y$ and boundary-changing line operator $\B$ from $\Y$ to a $\Z$ can be fused to produce a boundary-changing line operator from $\X$ to $\Z$. This operation is denoted $\otimes$.}
\end{figure}

In every set $\dW_{\X\X}$ there is a special element, the trivial boundary-changing line operator, which is an identity element with respect to the fusion product. The associativity of the fusion product is more difficult to formulate because there are boundary-changing line operators which are physically equivalent, but not equal. From the mathematical point of view, the difficulty can be explained  by saying that the set $\dW_{\X\Y}$ is actually a category, and it is not natural to talk about equality of objects in a category. Morphisms in this category are point operators inserted at points on the closed curve separating domains $\X$ and $\Y$ (Fig. 5). The insertion points of point operators  divide the closed curve into segments, and each segment is labeled by an element of the set $\dW_{\X\Y}$. 

\begin{figure}[htbp]  \label{fig5}
\centering
\includegraphics[bb=120 10 300 180,height=2in, scale=2.5, clip=true]{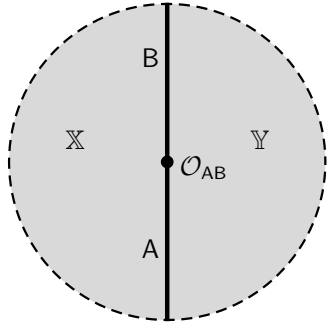}
\caption{Boundary-changing line operators between boundary conditions $\X$ and $\Y$ are objects of a category $\dW_{\X\Y}$. A morphism $\cO_{\A\B}$ from an object $\A$ to an object $\B$ is a point operator inserted at the junction of $\A$ and $\B$.}
\end{figure}

Let $\A$ and $\B$ be boundary-changing line operators between $\X$ and $\Y$; we will denote by $\dV_{\A\B}$ the set of point operators which can be inserted at the junction of segments labeled by $\A$ and $\B$. The same reasoning as in the case of  boundary-changing point operators tells us that $\dV_{\A\B}$ is a vector space. Point operators sitting on the same closed curve can be fused (Fig. 6), which gives rise to a product
$$
\dV_{\A\B}\times \dV_{\B\C}\ra {\dV}_{\A\C},
$$
which is bilinear and associative in an obvious sense. The category of boundary-changing line operators between $\X$ and $\Y$ has $\dW_{\X\Y}$ as its set of objects and the sets $\dV_{\A\B}$ as the sets of morphisms. Let us denote this category $\cC_{\X\Y}$. 

\begin{figure}[htbp]  \label{fig6}
\centering
\includegraphics[bb=100 10 550 200,height=2in, scale=2.5, clip=true]{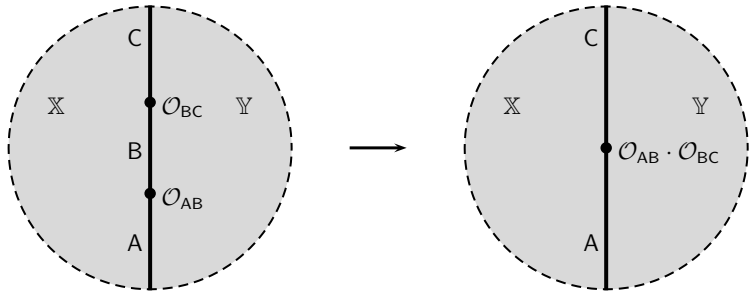}
\caption{Composition of morphisms in the category $\dW_{\X\Y}$ arises from the fusion of point operators sitting at the junctions of boundary-changing line operators.}
\end{figure}

The proper formulation of associativity of the fusion product (\ref{fusion}) says that that two triple products of three boundary-changing line operators differing by a placement of parentheses are isomorphic:
$$
(\A\otimes \B)\otimes\C \simeq \A\otimes (\B\otimes\C),\quad \forall \A\in\dW_{\X\Y},\forall \B\in\dW_{\Y\Z},\forall \C\in\dW_{\Z\T}.
$$
The isomorphism must be specified and must satisfy the so-called pentagon identity \cite{MacLane}. The above discussion can be summarized by saying that boundary conditions in a 3d TFT form a 2-category, whose sets of 1-morphisms are sets $\dW_{\X\Y}$ and whose sets of 2-morphisms are vector spaces $\dV_{\A\B}$. 

Kaluza-Klien reduction enables us to think about the category $\cC_{\X\Y}$ in two-dimensional terms. Consider a 3d TFT on a 3-manifold of the form $\Sigma\times [0,1]$ where 
$\Sigma$ is an oriented but not necessarily closed 2-manifold. On the boundaries $\Sigma\times \{0\}$ and $\Sigma\times \{1\}$ 
we impose boundary conditions $\X$ and $\Y$ respectively. Kaluza-Klein reduction tells us that one can describe the physics of this 3d TFT by an effective 2d TFT on 
$\Sigma$ which depends on $\X$ and $\Y$.  We claim that the category $\cC_{\X\Y}$ is the category of branes for this effective 2d TFT. Indeed, consider a 3d TFT on the half-space $\{(x,y,z)\vert x\geq 0\}$, where we imposed the boundary condition ${\mathbb X}$ on the half-plane 
$\{(0,y,z)\vert y<0\}$ and the boundary condition $\Y$ on the half-plane $\{(0,y,z)\vert y>0\}$. 
At the line given by $x=y=0$ we insert some boundary-changing line operator $\A\in \dW_{\X\Y}$. To regularize the problem we need to cut out a solid half-cylinder  $x^2+y^2<\eps^2$ for some $\eps>0$ and replace $\A$ with a suitable boundary condition on the part of the boundary given by $x^2+y^2=\eps^2$ (Fig. 7). Now we note that the half-space with a solid half-cylinder removed is diffeomorphic to $\RR_+\times\RR\times [-\pi/2,\pi/2]$, where $\RR_+$ is parameterized by the radial coordinate on the $(x,y)$ plane, $\RR$ is parameterized by $z$, and the interval is parameterized by the angular coordinate on the $(x,y)$ plane. Thus we may interpret the boundary condition at $x^2+y^2=\eps^2$ representing a boundary-changing line operator as a boundary condition in the effective 2d TFT on the half-space $\RR_+\times\RR$.

\begin{figure}[htbp]  \label{fig7}
\centering
\includegraphics[bb=40 10 390 180,height=2in, scale=2.5, clip=true]{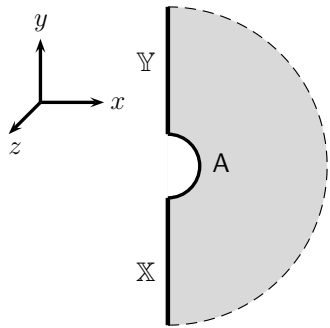}
\caption{A boundary-changing line operator is equivalent to a boundary condition on a half-cylinder $x^2+y^2=\eps^2,x>0$. This boundary condition is local in the $z$ direction but may be nonlocal in the angular direction in the $xy$ plane. It can be interpreted as a local boundary condition in a 2d TFT which is obtained by reducing the 3d TFT on an interval.}
\end{figure}

For any object $\X$ of a 2-category the endomorphism category $\cC_{\X\X}$ has a monoidal structure, i.e. an associative but not necessarily commutative tensor product. This monoidal structure is not natural from the 2d viewpoint (there is no physically reasonable way to define tensor product of branes in a general 2d TFT). Turning this around, if there is a mathematically natural monoidal structure on a category of branes in a 2d TFT, it is likely that this 2d TFT arises as a Kaluza-Klein reduction of a 3d TFT on an interval, and its category of branes can be interpreted as the category $\cC_{\X\X}$ for some boundary condition $\X$ in this 3d TFT. We will see an example of this below. 

As in 2d TFT, we may consider decorated manifolds, i.e. compact oriented manifolds with a nonempty boundary whose connected components are labeled by elements of the set of boundary conditions. Extended TFT in three dimensions assigns a number to a decorated 3-manifold (the value of the functional integral with given boundary conditions), a vector space to a decorated 2-manifold (the space of states of the effective 1d TFT obtained by Kaluza-Klein reduction on this 2-manifold), and a category to a decorated 1-manifold. The only connected decorated 1-manifold is an interval $[0,1]$.  If its endpoints are labeled by boundary conditions $\X$ and $\Y$, the corresponding category is the category of boundary-changing line operators $\cC_{\X\Y}$.

\subsection{Extended TFT in $n$ dimensions}

Continuing in the same fashion we conclude that boundary conditions in an $n$-dimensional TFT form an $n-1$-category. By analyzing more precisely the physical notion of a boundary condition one should be able to arrive at a physically-motivated definition of a weak $n$-category for all $n>0$. We will not try to do it here. 

We can now see why $n$-dimensional TFT assigns an $k$-category to a closed oriented $n-k-1$-manifold $M$. Consider an $n$-dimensional TFT on a manifold of the form $M\times N$, where $N$ is an oriented but not necessarily closed $k+1$-manifold. The Kaluza-Klein reduction principle tells us that we can describe the physics by an effective 
$k+1$-dimensional TFT on $N$. The $k$-category assigned to $M$ is the $k$-category of boundary conditions for this effective $k+1$-dimensional TFT.

If $M$ is the $n-k-1$-dimensional sphere, the corresponding $k$-category has an alternative interpretation: it is the category of defects of dimension $k$. To explain what a defect is, we might imagine that our TFT describes a particular macroscopic quantum state of a system of atoms in space-time of dimension $n$. It may happen that along some oriented submanifold $L$ of dimension $k$ the atoms are in a different state than elsewhere (or perhaps there is an altogether different kind of atoms inserted at this submanifold). In such a case one says that there is a defect of dimension $k$ inserted at $L$. Zero-dimensional defects are also known as local operators, one-dimensional defects are known as line operators, two-dimensional defects are known as surface operators. 

We claim that defects of dimension $k$ form a $k$-category, and that this category is nothing but the $k$-category assigned to $S^{n-k-1}$. To see this, suppose that $L$ is a $k$-dimensional plane in $\RR^n$. We introduce ``polar'' coordinates in $\RR^n$ such that $\RR^n\backslash L$ is identified with $\RR^k\times S^{n-k-1}\times\RR_+$, and $L$ is given by $r=0$, where $r$ is the coordinate on $\RR_+$. To regularize the problem we usually need to cut out a small tubular neighborhood of $L$ given by $r<\eps$ for some $\eps>0$ and replace the defect by a suitable boundary condition at $r=\eps$. Since our $n$-manifold has a factor $S^{n-k-1}$, we may regard this boundary condition as a boundary condition in an effective $k+1$-dimensional TFT which is obtained by Kaluza-Klien reduction on $S^{n-k-1}$. Thus defects of dimension $k$ can be regarded as objects of the $k$-category assigned to $S^{n-k-1}$. We will denote it $\cD_k$.

Let us note a few special cases. Local operators (i.e. defects of dimension $0$) are elements of the vector space assigned to $S^{n-1}$. This is usually called state-operator correspondence. Line operators are objects of a $\CC$-linear category, while surface operators are objects of a $\CC$-linear 2-category.

\begin{figure}[htbp]  \label{fig8}
\centering
\includegraphics[bb=120 10 300 200,height=2in, scale=2, clip=true]{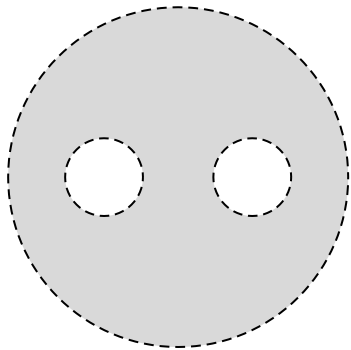}
\caption{A canonical bordism from $S^{n-k-1}\times S^{n-k-1}$ to $S^{n-k-1}$, here shown for $n-k=2$.}
\end{figure}

Two defects of dimension $k$ can be fused, which gives rise to a monoidal structure on $\cD_k$. Another way to understand the origin of this monoidal structure is to note that a solid ball of dimension $n-k$ with two smaller balls removed gives a canonical bordims from $S^{n-k-1}\times S^{n-k-1}$ and $S^{n-k-1}$ (Fig. 8). The Extended TFT assigns to this bordism a $k$-functor from $\cD_k\times\cD_k$ to $\cD_k$. This functor is associative (in some nontrivial sense). If $n-k$ is greater than $2$ it is also commutative, while for $n-k=2$ it is only braided. Thus monoidal, braided monoidal and symmetric monoidal $k$-categories also naturally arise from Extended TFT.  We note that M.~Kapranov and V.~Voevodsky proposed definitions for monoidal and braided monoidal 2-categories in \cite{KapVoe}. It is not clear if these definition agree with the definitions which are natural from the point of view of Extended TFT.

\section{The Rozansky-Witten model}

\subsection{Definition and basic properties}

The Rozansky-Witten model is a 3d TFT whose definition and basic properties have been described in \cite{RW}. It is a 3d sigma-model, i.e. its only bosonic field is a map $\phi:M\ra X$, where $M$ is an oriented 3-manifold, and $X$ is a complex manifold equipped with a holomorphic symplectic form. There are also fermionic fields: $\eta\in\Gamma(\phi^*T^{0.1}X)$ and $\rho\in\Gamma(\phi^*T^{1,0}X\otimes T^*M)$. The partition function of the RW model is defined as a functional integral over the infinite-dimensional supermanifold $\cM$ with local coordinates $(\phi,\eta,\rho)$. The measure is proportional to $\exp(-S(\phi,\eta,\rho))$, where $S$ is the action functional. Its explicit form is given in \cite{RW} but we will not need it here. The most important property of this measure is that it is invariant with respect to a certain odd vector field $Q$ on $\cM$ satisfying
$$
\{Q,Q\}=0,
$$
where braces denote the super-Lie bracket. As a result the partition function of the theory is unchanged if one adds to $S$ a function of the form $Q(f)$, where $f$ is any sufficiently well-behaved function on $\cM$. The Rozansky-Witten model is topological because its action can be written as
$$
S=Q(V)+S_0,
$$
where $S_0$ is independent of the metric on $M$ and $V$ is some metric-dependent odd function. 

The Rozansky-Witten model has a formal similarity to the better-known Chern-Simons gauge theory \cite{Shvarc,Witten:Jones}. For example, the non-$Q$-exact part of the action reads
$$
S_0=\int_M \left( \Omega(\rho,D\rho)+\frac{1}{3}\Omega(\rho,R(\rho,\rho,\eta))\right),
$$
where $\Omega$ is the holomorphic symplectic form on $X$, $D$ is the covariant differential on $\phi^*T^{1,0}X\otimes T^*M$ with respect to a pull-back of a connection on $T^{1,0}X$, and $R\in {\rm Hom}(T^{1,0}X\otimes T^{1,0}X\otimes T^{0,1}X,T^{1,0}X)$ is the curvature of this connection. Compare this to the action of the Chern-Simons gauge theory:
$$
S_{CS}=\int_M \left(\kappa(A,dA)+\frac{1}{3}\kappa(A,[A,A])\right),
$$
where $A$ is a gauge field (a connection on a principal $G$-bundle over $M$) and $\kappa$ is a non-degenerate $G$-invariant symmetric bilinear form on the Lie algebra of $G$. Clearly, the fermionic field $\rho$ should be regarded as analogous to the bosonic field $A$, i.e. the Rozansky-Witten model is an odd analogue of Chern-Simons theory. The symplectic form $\Omega$ is an analogue of the symmetric bilinear form $\kappa$.

This similarity is more than a mere analogy: as shown in \cite{RW} the Feynman diagram expansion of the partition function of the Rozansky-Witten model is basically the same as in Chern-Simons theory, with the curvature tensor playing the role of the structure constants of the Lie algebra of $G$. The partition function of the RW model (i.e. the number assigned to a closed oriented 3-manifold $M$) was shown in \cite{RW} to be a finite-type invariant of $M$. In the same paper the vector space assigned to a closed oriented 2-manifold of genus $g$ was computed; it turns out to be isomorphic to
$$
\bigoplus_p H^p_\bpartial \left(X,\left(\bigwedge T^{1,0}X\right)^{\otimes g}\right)
$$

The category associated to a circle turns out to be \cite{Roberts,RobWil} a 2-periodic version of the derived category of coherent sheaves on $X$ which we denote $\Dtwo(Coh(X))$. Its objects are 2-periodic complexes of coherent sheaves on $X$, and morphisms are obtained, as usual, by formally inverting quasi-isomorphisms. One way to see this is to consider the Rozansky-Witten model on $S^1\times\Sigma$, where $\Sigma$ is a not-necessarily-closed oriented 2-manifold. Kaluza-Klein reduction in this case gives a very simple effective 2d TFT: the so-called B-model with target $X$ which has been studied extensively in connection with Mirror Symmetry \cite{Witten:mirror}. Its category of branes is well understood and is known to be equivalent to the derived category of coherent sheaves \cite{Douglas,Mirrorbook}. The only difference compared to the usual B-model is that the $\ZZ$-grading is replaced by $\ZZ_2$-grading, leading to a 2-periodic version of the derived category.

The fact that the Rozansky-Witten model is $\ZZ_2$ -graded rather than $\ZZ$-graded might seem like a technicality, but it becomes crucially important when one turns to computing the braided monoidal structure on the category assigned to $S^1$. Both the usual and the 2-periodic derived categories have obvious symmetric monoidal structures given by the derived tensor product. However the braided monoidal structure on $\Dtwo(Coh(X))$ which arises from the Extended TFT turns out to be not this obvious symmetric monoidal structure but its deformation, which is no longer symmetric \cite{Roberts,RobWil}. $\ZZ_2$-grading is important here, because the usual bounded derived category appears not to admit any non-symmetric monoidal deformations (see below). 

The braided monoidal deformation of $\Dtwo(Coh(X))$ is a quantum deformation of the obvious monoidal structure, in the same sense that the category of representations of the quantum group is a quantum deformation of the category of representations of the corresponding classical group. That is, corrections to the ``obvious'' associator arise as quantum corrections in the Feynman diagram expansion of the Rozansky-Witten model. This has been worked out in detail in \cite{RobWil}. The analogue of the Planck constant is the inverse of the symplectic form $\Omega$, i.e. under a rescaling $\Omega\ra \lambda\Omega$ the $p$-th order quantum correction scales like $\lambda^{1-p}$. 

Finally let us turn to the 2-category that the Rozansky-Witten model assigns  to a point, i.e. the 2-category of boundary conditions \cite{KRS,KR}. Its simplest objects are complex Lagrangian submanifolds of $X$. There are also more complicated objects which are described by a family of Calabi-Yau manifolds parameterized by points of a complex Lagrangian submanifold $Y$, or even more abstractly, a family of Calabi-Yau categories (i.e. categories which have the same formal properties as the derived category of coherent sheaves of a Calabi-Yau manifold) over $Y$. For simplicity we will not discuss these more complicated objects here. 

Even for the simplest geometric objects the description of the categories of morphisms turns out quite complicated, so we will discuss only two extreme cases (see \cite{KR} for a more general case). First, suppose that two complex Lagrangian submanifolds $Y_1$ and $Y_2$ intersect at isolated points, but not necessarily transversely. The category of morphisms in this case is the direct sum of categories corresponding to each intersection point, in the sense that the set of objects is the union of the sets of objects assigned to each point, and spaces of morphisms between objects coming from different intersection points are 0-dimensional vector spaces. Thus it suffices to describe the category corresponding to a single intersection point. 

In the neighborhood of the intersection point one may choose complex Darboux coordinates, i.e. choose an identification of the neighborhood with an open subset $U$ of $T^*\CC^m$ with its canonical symplectic form. One can always choose this identification so that $Y_1\cap U$ and $Y_2\cap U$ are represented by graphs of exact holomorphic 1-forms $dW_1$ and $dW_2$ where $W_1$ and $W_2$ are holomoprhic functions on an open subset $V$ of $\CC^m$. The category of morphisms in this case is equivalent to the {\it category of matrix factorizations} of the function $W_2-W_1$ \cite{KRS,KR}. This category was introduced by M.~Kontsevich in connection with Homological Mirror Symmetry and is defined as follows. An object  of this category is a $\ZZ_2$-graded holomorphic vector bundle $E$ on $V$ equipped with a holomorphic endomorphism $D$ of odd degree satisfying $D^2=W_2-W_1+c$, where $c$ is a complex number. Consider two such objects $(E_1,D_1,c_1)$ and $(E_2,D_2,c_2)$ The space of holomorphic bundle maps  from $E_1$ to $E_2$ is a $\ZZ_2$-graded vector space equipped with an odd endomorphism $D_{12}$ defined by
$$
D_{12}(\phi)=D_2\cdot \phi - (-1)^{|\phi|} \phi\cdot  D_1,\quad \forall \phi\in \Hom(E_1,E_2),
$$
where $|\phi|=0$ or $1$ depending on whether $\phi$ is even or odd. The endomorphism $D_{12}$ satisfies $D_{12}^2=c_2-c_1$. The space of morphisms in the category of matrix factorizations from $(E_1,D_1,c_1)$ to $(E_2,D_2,c_2)$ is defined to be the cohomology of the differential $D_{12}$ if $c_1=c_2$; otherwise it is defined to be zero. 

The other relatively simple case is the category of endomorphisms of a complex Lagrangian submanifold $Y$ of $X$. On the classical level it is equivalent to $\Dtwo(Coh(Y))$. One way to see this is to perform the Kaluza-Klein reduction of the Rozansky-Witten model on an interval $[0,1]$ with boundary conditions corresponding to $Y$. On the classical level the effective 2d TFT is a $\ZZ_2$-graded version of the B-model with target $Y$, whose category of branes is $\Dtwo(Coh(Y))$. The monoidal structure is the usual derived tensor product. Both the monoidal structure and the category itself are modified by quantum corrections, in general.

\subsection{Monoidal deformations of the derived category of coherent sheaves}

As far as we know, the theory of deformations of monoidal categories in general, and the derived category of coherent sheaves in particular, has not been systematically developed. A remarkably simple geometric picture of monoidal deformations of the latter category emerges from the study of the Rozansky-Witten model. Consider a complex Lagrangian submanifold $Y$ of a holomorphic symplectic manifold $(X,\Omega)$. The functional integral of the Rozansky-Witten model localizes on constant maps, which implies that the category of endomorphisms of the object $Y$ depends only on the formal neighborhood of $Y$ in $X$. If this formal neighborhood happens to be isomorphic, as a holomorphic symplectic manifold, to the formal neighborhood of the zero section of $T^*Y$ (which we will denote $T^*_fY$), then one can show that the category of endomorphisms of $Y$ does not receive quantum corrections and therefore is equivalent to $\Dtwo(Coh(Y))$ as a monoidal category. 

In general, the formal neighborhood of $Y$ is isomorphic to $T^*_fY$ as a real symplectic manifold, but not as a holomorphic symplectic one.  The deviation of the complex structure on $T^*_fY$ from the standard one is described by a $(0,1)$-form $\beta$ with values in the graded holomorphic vector bundle 
\begin{equation}\label{symp2}
\oplus_{p=2}^\infty  \Sym^p(TY).
\end{equation}
This $(0,1)$-form satisfies a Maurer-Cartan-type equation
\begin{equation}\label{MC}
\bpartial\beta+\frac12 [\beta,\beta]=0.
\end{equation}
Here brackets denote wedge product of forms combined with a Lie bracket on sections of $\Sym^\bul TY$. The Lie bracket comes from the identification of the space of sections of $\Sym^\bul TY$ with the space of fiberwise-holomorphic functions on $T^*_f Y$ and the Poisson bracket on such functions. Note that because the wedge product of 1-forms is skew-symmetric, the expression $[\beta,\beta']$ is symmetric with respect to the exchange of 1-forms $\beta$ and $\beta'$, and therefore $[\beta,\beta]$ need not vanish.

The Rozansky-Witten model provides a map from the space of solutions of the equation (\ref{MC}) to the space of monoidal deformations of $\Dtwo(Coh(Y))$. As usual, there is a group of gauge transformations whose action on the space of solutions is determined by the action of its Lie algebra:
\begin{equation}\label{gaugeMC}
a: \beta\mapsto \beta+\bpartial a+[\beta,a],
\end{equation}
where $a$ is a section of $\oplus_{p=1}^\infty \Sym^p TY$. 

The following natural conjecture was formulated in \cite{KRS}:
\begin{conjecture}
Let $\M_Y$ be the space of solutions of the Maurer-Cartan equation (\ref{MC}) where $\beta$ is an inhomogeneous form of type $(0,q)$ on $Y$ with odd $q$ with values in the bundle (\ref{symp2}). There is a surjective map from $\M_Y$ to the space of monoidal deformations of the category $\Dtwo(Coh(Y))$. Two monoidal deformations are equivalent if the corresponding solutions of the Maurer-Cartan equation are related by a gauge transformation whose infinitesimal form is given by (\ref{gaugeMC}).
\end{conjecture} 
If true, this conjecture gives an elegant description of all monoidal deformations of $\Dtwo(Coh(Y))$.  This description is strikingly similar to the description of all deformations, monoidal or not, of the category $D^b(Coh(Y))$ regarded as an $A_\infty$ category \cite{Konts:Formality}. The latter makes use of a Maurer-Cartan-type equation for a $(0,q)$ form $P$ with values in the graded bundle $\Lambda^\bul TY$:
\begin{equation}\label{BKMC}
\bpartial P+\frac12 [P,P]_{SN}=0. 
\end{equation}
Here brackets denote the Schouten-Nijenjuis bracket on polyvector fields, and the total degree of $P$ (that is, the sum of the form degree and the polyvector degree) is even. In particular a holomorphic Poisson bivector, i.e. a section of $\Lambda^2 TY$ satisfying
$$
\bpartial P=0,\quad [P,P]_{SN}=0
$$
gives rise to a noncommutative deformation of $Y$. The analog of the holomorphic Poisson bivector in our case is a$(0,q)$-form $\beta$ with values in $\Sym^2 TY$ satisfying
$$
\bpartial\beta=0,\quad [\beta,\beta]=0.
$$
The corresponding deformation of the monoidal structure makes the tensor product on $\Dtwo(Coh(Y))$ non-symmetric \cite{KRS,KR}. Thus one may regard this deformation as a categorification of deformation quantization, and the conjectural relation between the space of solutions of the Maurer-Cartan equation (\ref{MC}) and the space of monoidal deformations as a categorification of the Formality Theorem of M.~Kontsevich \cite{Konts:Formality}.

Let us comment on the analogue of the above conjecture in the $\ZZ$-graded case, i.e. when $\Dtwo(Coh(Y))$ is replaced with $D^b(Coh(Y))$. The latter category can be interpreted as the endomorphism category of a boundary condition in a $\ZZ$-graded version of the Rozansky-Witten model. This $\ZZ$-graded version exists if the target manifold $X$ admits a $\CC^*$ action with respect to which the holomorphic symplectic form has weight $2$. To realize $D^b(Coh(Y))$ as the endomorphism category, we take $X=T^*Y$ with a canonical symplectic form $dp dq$, where $p$ is the fiber coordinate, and define the $\CC^*$ action by
$$
\lambda: p\mapsto \lambda^2 p,\quad \lambda\in\CC^*.
$$
That is, the fiber coordinate has weight $2$. Accordingly, if we identify the space of sections of $\Sym^p TY$ with the space of functions on $T^*Y$ which are holomorphic degree-$p$ polynomials on the fibers, we should place it in cohomological degree $2p$. A $(0,q)$-form with values in $\Sym^p TY$ will therefore have degree $q+2p$. From the point of view of the $\ZZ$-graded Rozansky-Witten model, holomorphic symplectic deformations of $T^*Y$ should be identified with degree-$3$ solutions of the Maurer-Cartan equation (\ref{MC}). Since $p\geq 2$, the only such solution is the zero one.\footnote{The restriction to $p\geq 2$ appears because we want the complex structure of $Y$ to be undeformed. If we relax this assumption, then the only allowed $\beta$ is a $(0,1)$-form with values in $TY$ and satisfying the Maurer-Cartan equation (\ref{BKMC}). Such $\beta$ describes a deformation of $D^b(Coh(Y))$ which arises from a deformation of the complex structure on $Y$. This deformation is obviously monoidal, but not very interesting.}

The $\ZZ$-graded Rozansky-Witten model with target $T^*Y$ is particularly simple since one can show that quantum corrections always vanish. This enables one to give a fairly concise description of the 2-category of boundary conditions for this 3d TFT. If we regard a monoidal category $\cC$ as a categorification of an associative algebra, then the categorification of a module is a module category over $\cC$, .i.e.  a category $\cD$ on which $\cC$ acts by endofunctors. Module categories over a monoidal category form a 2-category. Roughly speaking, the 2-category of boundary conditions for the $\ZZ$-graded Rozansky-Witten model with target $T^*Y$ is the 2-category of module categories over the monoidal category $D^b(Coh(Y))$. To be more precise, one needs to consider a differential graded version of both the monoidal category $D^b(Coh(Y))$ and module categories over it. The resulting 2-category also appeared in the mathematical papers \cite{BZFN,ToVe}. 

Note that the existence of the $\ZZ$-graded Rozansky-Witten model may be regarded as a ``physical reason'' for the existence of the derived tensor product on the category $D^b(Coh(Y))$ associated to a complex manifold $Y$. That is, while the derived tensor product has no physical meaning if one thinks about $D^b(Coh(Y))$ as the category of branes in a 2d TFT (the B-model with target $Y$), it arises naturally once we realize that the B-model with target $Y$ can be obtained by Kaluza-Klein reduction from a 3d TFT on an interval, with identical boundary conditions at the two endpoints. This 3d TFT is the $\ZZ$-graded Rozansky-Witten model with target $T^*Y$, and the boundary condition is the complex Lagrangian submanifold given by the zero section of $T^*Y$. 

\section{Topological Gauge Theory in four dimensions and Geometric Langlands Duality}

\subsection{Electric-magnetic duality and Topological Gauge Theory}

Another application of the Extended TFT has to do with the Geometric Langlands Duality. \footnote{The relationship between ETFT and Geometric Langlands Duality is also discussed in \cite{BZN}.} The physical origin of the Geometric Langlands Duality is a conjectural isomorphism between two supersymmetric gauge theories in four dimensions with gauge groups $G$ and $\LG$, where $\LG$ is the Langlands-dual of $G$.\footnote{We remind that the Langlands dual of compact simple Lie group $G$ is a compact simple Lie group $\LG$ whose maximal torus is isomorphic to the dual of the maximal torus of $G$.} This isomorphism is known as Montonen-Olive duality \cite{MO} and holds for gauge theories with maximal supersymmetry \cite{WittenOlive, Osborn}. There are many computations verifying particular implications of the Montonen-Olive conjecture, but no general proof. In the case $G=U(1)=\LG$, the Montonen-Olive duality reduces to electric-magnetic duality, i.e. the transformation which exchanges electric and magnetic fields. Electric-magnetic duality is well known to be a symmetry of the $U(1)$ gauge theory both on the classical and quantum levels. The Montonen-Olive conjecture can be regarded as a far-reaching nonabelian generalization of this  fact. 

To connect the Montonen-Olive conjecture to Geometric Langlands Duality the first step is to replace supersymmetric gauge theories with much simpler topological field theories \cite{KW}. This is achieved by means of a procedure called {\it twisting} \cite{Witten:DonaldsonTFT}. Roughly speaking, one redefines the stress-energy tensor of the theory, so that it becomes $Q$-exact with respect to a certain nilpotent odd vector field $Q$ on the supermanifold of all field configurations. For historical reasons, $Q$ is known as the BRST operator. Simultaneously one restricts observables (i.e. functions on the supermanifold of field configurations) to those which are annihilated by $Q$. After twisting both the supersymmetric gauge theory with gauge group $G$ and its Montonen-Olive dual with gauge group $\LG$ one gets a pair  of isomorphic 4d TFTs. 

In fact, the situation is more complicated than that. First of all, the maximally supersymmetric gauge theory in four dimensions admits three inequivalent twists differing both by the choice of $Q$ and by the required modification of the stress-energy tensor \cite{VafaWitten}. Following \cite{KW} we will focus on the twist which was first considered by N.~Marcus \cite{Marcus} and is nowadays called the GL-twist. Second, the GL-twisted gauge theory has two supercommuting BRST operators $Q_l$ and $Q_r$, and one can take any linear combination of this as the BRST operator which must annihilate the observables:
$$
Q=u Q_l+v Q_r,\quad u,v \in \CC,\quad |u|^2+|v|^2>0. 
$$
The overall normalization of $Q$ does not affect the theory, so the GL-twisted theory is really a family of 4d TFTs parameterized by points of $\CC\PP^1$ \cite{KW}. We will identify $\CC\PP^1$ with a one-point compactification of the complex plane and will label a particular TFT by a parameter $t\in \CC\cup \{\infty\}$.

It turns out that Montonen-Olive duality not only replaces $G$ with $\LG$ but also acts on the parameter $t$ \cite{KW}. Geometric Langlands duality arises from a particular instance of the Montonen-Olive duality which maps $t=i$ to $t=1$.\footnote{More generally, one gets what is known as Quantum Geometric Langlands.}

\subsection{From Topological Gauge Theory to Geometric Langlands Duality}

We can now attempt to extract some mathematical consequences of the Montonen-Olive conjecture. Let us fix a compact simple Lie group $G$ and let $C$ be a closed oriented 2-manifold.  The 4d Topological Gauge Theory with gauge group $G$ assigns to $C$ a family of $\CC$-linear categories parameterized by $t$; we will denote a member of this family  by $\F(G,t,C)$. The Montonen-Olive conjecture implies that there is an equivalence of categories
$$
\F(G,i,C)\simeq \F(\LG,1,C).
$$
It remains to understand the categories involved. This turns out to be rather nontrivial. The category $\F(G,t,C)$ is the category of branes for the 2d TFT obtained by Kaluza-Klein reduction of the 4d Topological Gauge Theory on $C$. For $g>1$ this 2d TFT was analyzed in \cite{KW}. It was shown that for $t=i$ the 2d TFT is a B-model whose target is the moduli space $\cM_{flat}(G_\CC,C)$ of flat $G_\CC$-connections on $C$.  Here $G_\CC$ is the complexification of $G$. Accordingly, for $t=i$ the category of branes of the 2d TFT is the derived category of coherent sheaves on $\cM_{flat}(G_\CC,C)$. For $t=1$ the 2d TFT is a different topological sigma-model (the A-model) whose target space is a symplectic manifold $\cM_{flat}^{symp}(G_\CC,C)$. This manifold is diffeomorphic to $\cM_{flat}(G_\CC,C)$, with an exact symplectic form given by
$$
\omega=\int_C \kappa(\delta A\wedge \delta\phi),
$$
where $A$ is the real part of the flat $G_\CC$-connection on $C$, $\phi$ is its imaginary part, and $\kappa$ is the Killing metric on the Lie algebra of $G$. The category of branes for an A-model is the so-called Fukaya-Floer category \cite{FOOO} whose objects are Lagrangian submanifolds of the target space and morphisms are defined by means of the Lagrangian Floer homology. Thus Montonen-Olive duality implies that the derived category of coherent sheaves on $\cM_{flat}(G_\CC,C)$ is equivalent to the Fukaya-Floer category of $\cM_{flat}^{symp}(\LG_\CC,C)$. 

The usual statement of the Geometric Langlands Duality is somewhat different. Instead of the Fukaya-Floer category of $\cM_{flat}^{symp}(\LG_\CC,C)$ it involves the derived category of D-modules over the moduli stack of holomorphic $\LG_\CC$-bundles over $C$. It was shown in \cite{KW} that there is a functor from the former to the latter, but is not clear at the time of writing why this functor should be an equivalence. 

If $C$ has genus zero, the effective 2d TFT one obtains by Kaluza-Klein reduction is rather different. For $t=i$ it was shown in \cite{KSV} that the 2d TFT is a $G_\CC$-equivariant B-model whose target is the Lie algebra $\frg_\CC$ of $G_\CC$ placed in cohomological degree $2$. That is, the target is a purely even graded manifold which we denote $\frg_\CC[2]$. The group $G_\CC$ acts on $\frg_\CC$ by the adjoint representation. The category of branes for this 2d TFT is the $G_\CC$-equivariant derived category of coherent sheaves on $\frg_\CC[2]$ \cite{KSV}. For $t=0$ and genus zero the 2d TFT has not been analyzed thoroughly yet. From the mathematical viewpoint, one expects that Geometric Langlands Duality relates the category $D^b_{G_\CC}(Coh(\frg[2]))$ and the $\LG$-equivariant constructible derived category of sheaves on the loop group of $\LG$ \cite{BezFink}. It would be interesting to see if the latter category emerges as the category of branes in a 2d TFT. 

Instead of studying categories which the 4d TFT assigns to a closed oriented 2-manifold, we may consider 2-categories assigned to a circle. For a fixed gauge group $G$ we get a family of these 2-categories parameterized by $t$. Let us denote a member of this family $\F(G,t,S^1)$. The Montonen-Olive conjecture implies an equivalence of 2-categories:
\begin{equation}\label{GL:twocat}
\F(G,i,S^1)\simeq \F(\LG,1,S^1).
\end{equation}
Moreover, as mentioned above, both 2-categories are supposed to have braided monoidal structure, and the equivalence is supposed to be compatible with them. 

To understand the 2-categories $\F(G,t,S^1)$ one needs to study the Kaluza-Klein reduction of the 4d Topological Gauge Theory on a circle. The corresponding 3d TFT has been analyzed in \cite{KSV}. It turns out that at $t=i$ the 3d TFT is a $G_\CC$-equivariant version of the $\ZZ$-graded Rozansky-Witten model with target $T^*G_\CC$, where $G_\CC$ acts on $T^*G_\CC$ by conjugation. The corresponding 2-category of branes is, roughly speaking, the 2-category of module categories over the $G_\CC$-equivariant derived category of coherent sheaves over $G_\CC$, regarded as a monoidal category. A typical object of this 2-category is a family of categories over $G_\CC$ with an action of $G_\CC$ which lifts the action of $G_\CC$ on itself by conjugation. For $t=1$ the 3d TFT has not been studied thoroughly, but it appears that its objects are module categories over a $G$-equivariant version of the Fukaya-Floer category of $T^*G$, regarded as a monoidal category.\footnote{The monoidal structure on this category is less obvious than for the analogous category at $t=i$.}

\section{Open questions}

It  is plausible that the correct mathematical formalism for Extended Topological Field Theories is provided by the theory of $(\infty,n)$-categories developed by J.~Lurie \cite{Lurie:TFT,Lurie:Goodwillie}. More precisely, what seems most relevant is the linear case of this theory where the set of $n$-morphisms has the structure of a differential graded vector space.\footnote{From the physical viewpoint, linearity arises from the fact that the quantum-mechanical space of states is a (graded) vector space. } Specifically, one expects that to every $n$-dimensional TFT one can attach a linear $(\infty,n-1)$-category whose objects are boundary conditions, whose 1-morphisms are boundary-changing operators supported on submanifolds of the boundary of codimension $1$, etc.  Compositions of $k$-morphisms arises from the fusion of the physical operators supported on submanifolds of the boundary of codimension $k$. It would be interesting to understand whether the definition of the $(\infty,n)$-category captures the properties of fusion expected on physical grounds. 

Extended Topological Field Theory provides a new viewpoint on the Geometric Langlands Program. The most powerful statement implied by the topologically twisted version of the Montonen-Olive conjecture is the equivalence of 3-categories which the 4d Topological Gauge Theories with gauge groups $G$ and $\LG$ assign to a point. From the physical viewpoint, these are 3-categories of boundary conditions. Some examples of boundary conditions for maximally supersymmetric gauge theories and their Montonen-Olive duals have been constructed in \cite{GaWi1,GaWi2}, but we are still far from understanding the nature of this 3-category. 

One categorical level down, it would be very interesting to study the rich structure on the 2-category $\F(G,t,S^1)$ which the 4d Topological Gauge Theory assigns to $S^1$.  We already mentioned that it has a braided monoidal structure and an identity object. The monoidal structure arises from the bordism from $S^1\times S^1$ to $S^1$ depicted in fig. 8. Another way to draw it is shown in fig. 9; we will call it the ``pants'' bordism. Similarly, the identity object arises from a disc regarded as a bordism between the empty 1-manifold and $S^1$ (the ``cup'' bordism, see fig. 9). Further, the 2-category $\F(G,t,S^1)$ is expected to be {\it rigid}, i.e. there should be a 2-functor $e$ from $\F(G,t,S^1)\times \F(G,t,S^1)$ to the 2-category of linear categories, and a 2-functor $\iota$ in the opposite direction satisfying some compatibility conditions. These 2-functors both arise from a cylinder regarded either as a bordism from $S^1\times S^1$ to the empty 1-manifold (the "downward plumbing fixture", see fig. 9), or as a bordism from the empty 1-manifold to $S^1\times S^1$ (the "upward plumbing fixture", see fig. 9). All these 2-functors should satisfy various compatibility conditions arising from the fact that although one can glue a given oriented 2-manifold with boundaries from these four building blocks ("pants", "cup" and two "plumbing fixtures") in many different ways, the equivalence class of 2-functors corresponding to this 2-manifold must be well-defined.  One can summarize the situation by saying that $\F(G,t,S^1)$ is a rigid braided monoidal 2-category.  A related viewpoint on the Geometric Langlands Duality (not using the language of Extended Topological Field Theory) is proposed in \cite{FG}.

\begin{figure}[htbp]  \label{fig9}
\centering
\includegraphics[bb=100 30 550 125,height=1.1in, scale=2.5, clip=true]{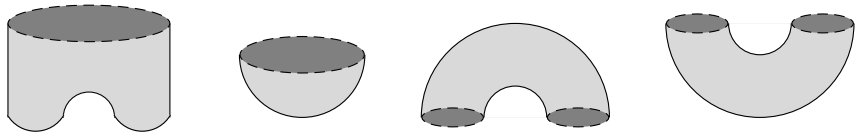}
\caption{Basic 2d bordisms: pants, cup, downward plumbing fixture, upward plumbing fixture.}
\end{figure}

Montonen-Olive duality predicts that $\F(G,i,S^1)$ and $\F(\LG,1,S^1)$ are equivalent as rigid braided monoidal 2-categories. This statement should imply the statement of the usual Geometric Langlands Duality in the following way. Given any closed oriented 2-manifold $C$ we may represent it as a result of gluing the four building blocks shown in fig. 9. Extended TFT in four dimensions assigns to every building block a 2-functor as described above, and therefore assigns to the whole $C$ a 2-functor from the 2-category of linear categories to itself. Axioms of Extended TFT ensure that the equivalence class of this 2-functor is independent of the way one cut $C$ into pieces. The category $\F(G,t,C)$ can be thought of as the result of applying this 2-functor to the category of vector spaces. The rigid braided monoidal equivalence of $\F(G,i,S^1)$ and $\F(\LG,1,S^1)$ then implies the equivalence of categories $\F(G,i,C)$ and $\F(\LG,1,C)$.

\section*{Acknowledgments} 
I would like to thank David Ben-Zvi, Dmitri Orlov, and Lev Rozansky for comments on the draft of the talk.

\end{document}